\documentstyle[twoside]{article} 
\newtheorem{thm}{Theorem}
\newtheorem{example}{Example}
\newtheorem{defn}[thm]{Definition}
\newtheorem{lemma}[thm]{Lemma}
\newtheorem{prop}[thm]{Proposition}
\newtheorem{cor}[thm]{Corollary}

\newcommand{\fl}[1]{\lfloor #1 \rfloor}
\newcommand{\nin}{\not\in}
\begin{document} 

\begin{titlepage} 
\title{Decompositions of ${\cal B}_{n}$ and $\Pi _{n}$ using Symmetric
Chains} 
\author{D. Loeb\\
D\'{e}partement de Math\'{e}matiques et d'Informatique\\
Universit\'{e} de Bordeaux I\\
33405 Talence, France \and 
E. Damiani and O. D'Antona\\
Dipartimento di Scienze dell'Informazione\\
Universit\`a di Milano\\
Milano, Italy}
\end{titlepage} 
\maketitle

\begin{abstract}
We review the Green/Kleitman/Leeb interpretation of de Bruijn's
symmetric chain decomposition of ${\cal B}_{n}$, and explain how it
can be used to find a maximal  
collection of disjoint symmetric chains in the nonsymmetric lattice of
partitions of a set.
\end{abstract}

\section{Introduction}

De Bruijn {\em et al} \cite{DeB} and Griggs \cite{Griggs} have given
symmetric chain 
decompositions of various lattices including in particular the Boolean
lattice  of subsets of a given set. 
Simpler characterizations of this decomposition were later found by
Greene and Kleitman \cite{GK} and independently by K. Leeb
[unpublished]. These results are of great utility in combinatorics,
permitting, for example, various extensions of Sperner's theorem \cite{S}.

After a review of previous work for the Boolean lattice, we
investigate whether a similar decomposition can be found for the
lattice of partitions of a set.
This lattice is not
symmetric. For example, there are $2^{n-1}-1$ partitions of an $n$-set
into two blocks, but only $n(n-1)/2$ partitions of an $n$-set into
$n-1$ blocks. Thus, there is no complete decomposition of the lattice
into symmetric chains (for $n\geq 3$). However, using de Bruijn's
decomposition together with a certain method to encode sets,
we find a maximal collection of disjoint symmetric chains of
partitions.  

We conclude with other applications of this set encoding: it gives the
relationship between the elementary and complete symmetric functions
as well as a new formula for the Bell numbers.

\section{de Bruijn's Method}

A {\em symmetric chain} in a ranked lattice $L$
is a sequence of elements $(x_{i})_{i=1}^{n}$ such that $x_{i+1}$
covers $x_{i}$ and $r(x_{1})+r(x_{n})=r(L)$. 
A {\em partition} $\phi $ of a set $S$ is an unordered collection of disjoint
subsets called {\em blocks} whose union is $S$. (We write $\phi \vdash
S.$) A {\em decomposition} of
a lattice $L$ into symmetric chains is a partition of $L$ whose blocks are
symmetric chains.

\begin{prop}[de Bruijn \cite{DeB}]
Let $L_{1}$ and $L_{2}$ be lattices with symmetric chain
decompositions. Then $L_{1}\times L_{2}$ has a symmetric chain
decomposition.$\Box $
\end{prop}

{\em Proof:} Let ${\cal D}_{1}$ and ${\cal D}_{2}$ be decompositions
of $L_{1}$ and $L_{2}$ respectively. It suffices to decompose $C\times
D$ into symmetric chains for every $C\in {\cal D}_{1}$ and
$D\in {\cal D}_{2}$. Suppose that $C = \{c_{i} \}_{i=1}^{k}$ and
$D=\{d_{j} \}_{j=1}^{l}$. Then the sets $\{E_{j} \}_{j=0}^{l}$ are such
a symmetric chain decomposition, where the individual chains $E_{j}$
are given by 
$$ E_{j} = \{(c_{1},d_{j}) < \ldots < (c_{k-j},d_{j}) <
(c_{k-j},d_{j+1}) < \ldots < (c_{k-j},d_{l})  \}.\Box  $$

\begin{cor}[de Bruijn \cite{DeB}]
Any finite chain product is symmetric (e.g.: ${\cal B}_{n}$).$\Box$
\end{cor}

Let ${\cal D}_{n}$ be the decomposition of ${\cal B}_{n}$ given by
de Bruijn. Then ${\cal D}_{n+1}$ can be constructed from ${\cal D}_{n}$
as follows: 
$${\cal D}_{n+1}= \alpha({\cal D}_{n}) \cup \beta ({\cal D}_{n}) -
\{\emptyset  \}   $$
where  $\emptyset $ denotes the empty chain and $\alpha$ and $\beta $
are the maps on chains: 
\begin{eqnarray*}
\alpha( c_{1},\ldots,c_{k}) &=& c_{1}\cdots c_{k}(c_{k}\cup \{n+1 \})\\
\beta ( c_{1},\ldots,c_{k}) &=& (c_{1}\cup \{n+1 \})\cdots (c_{k}\cup
\{n+1 \}).
\end{eqnarray*}

Greene, Kleitman, Leeb give \cite{GK} a direct description of these symmetric
chains. This description is equivalent to de Bruijn's; however, it has
the advantage that a direct comparison of two sets can indicate
whether they are in the same chain, and a complete chain can easily be
generated given one of its elements.

We write $S\subseteq \{1,2,\ldots,n \}$ as a word $w(S)$ of length
$n$ with a right or left parenthesis in position $i$  according to
whether $i\in S$,
$$ w(S)_{i}= \left\{ \begin{array}{ll}
)		&\mbox{if $i\in S$}\\
(		&\mbox{if $i\nin S$.}
         \end{array}
\right. $$
For example, for $n=10$ and  
 $S=\{1,3,4,8,9 \}$,
$w(S)$ is )())((())(. We can then indicate matching parentheses in
bold: $w(S)=$ ){\bf{}()})({\bf{}(())}(. Those sets with the same
matching parentheses belong 
to the same chain. Moreover, it is easy to see how to generate the rest of
the chain. The parentheses in roman face must not match. Thus, they
must consist of a certain number of right parentheses followed by a
certain number of left parentheses. Thus, 
$$ \begin{array}{ll}
\mbox{({\bf{}()}(({\bf{}(())}(},\\
\mbox{){\bf{}()}(({\bf{}(())}(},\\
\mbox{){\bf{}()})({\bf{}(())}(},\\
\mbox{){\bf{}()})){\bf{}(())}(},&\mbox{and}\\
\mbox{){\bf{}()})){\bf{}(())})}
\end{array}$$
correspond respectively to the sets 
$$\begin{array}{ll}
\{3,8,9 \},\\
\{1,3,8,9 \},\\
\{1,3,4,8,9 \}=S,\\
\{1,3,4,5,8,9 \},&\mbox{and}\\
\{1,3,4,5,8,9,10\}.
\end{array}$$

We note several easily verified facts.
\begin{prop}\label{starprop}
\begin{enumerate}
\item In each chain, elements are added in increasing order.
\item $n$ is an element of the last set in each chain.
\item If $S\rightarrow S\cup \{i \}$ is a link in a chain, then
$i+1\nin S$, and either $i=1$ or $i-1\in S.\Box $
 \label{star}
\end{enumerate}
\end{prop}

One should note, however, that not all symmetric chain decompositions
are equivalent, even up to automorphism. For example,
consider the De Bruijn decomposition of ${\cal B}_{4}$. The pair of
chains: 
$$(\{4 \},\{1,4 \},\{1,2,4\})\mbox{ and }(\{2,4 \})$$
can be replaced by the pair of chains  
$$(\{4 \},\{2,4 \},\{1,2,4\})\mbox{ and }(\{1,4 \}).$$
The resulting decomposition is not isomorphic to the original.
 
See \cite{MR,MRSW} for other systems of ``parentheses'' used to
find a Sperner decomposition of a poset into chains.

\section{Partition Lattice}

We will now describe how de Bruijn's  decomposition of
${\cal B}_{n}$ into symmetric chains can be used to generate a similar
decomposition of the lattice $\Pi _{n+1}$
of partitions of the set $\{1,2,\ldots,n+1 \}$.

 A partition $\phi \vdash S$ is
a {\em refinement} of another partition $\kappa \vdash S $ if each block of
$\kappa $ is a union of blocks of $\phi $. This ordering results in a
ranked lattice $\Pi _{|S|}$.
Partitions $\pi \in \Pi _{n}$ of rank $i$ consist of $n-i$ blocks. 

Since there are $n(n-1)/2$ partitions with $n-1$ blocks and
$2^{n-1}-1$ partitions with two blocks, there is no total
decomposition of $\Pi _{n}$ into symmetric chains. However, amazingly
enough we can find a set of disjoint symmetric chains which includes
all partitions of rank $\leq \fl{(n-1)/2}$. Such a collection is
clearly maximal.

At first glance, such a collection may be surprising. However, it is
not a Sperner 
decomposition, so it does not contradict the fact that $\Pi
_{n}$ is neither LYM \cite{JHS} nor Sperner \cite{ERC,JBS,JK}.
\begin{defn}[Coding of Sets] \label{coding}
For each $S\subseteq \{1,2,\ldots,n \}$, we associate a code $c(S)\in
{\bf N}_{0}^{n+1}$ of length $n+1$ as follows:
$$c(S)_i = \left\{ \begin{array}{ll}
0		&\mbox{if $i\in S$}\\
i-\sum _{j=1}^{i-1}c(S)_{j}		&\mbox{if $i\nin S$.}
         \end{array}\right.$$
\end{defn}

For example, the code of $S =\{1,2,3,7,11,12,16,18,19\}\subseteq
\{1,2,\ldots,20 \} $ is
$$c(S)=(0,0,0,4,1,1,0,2,1,1,0,0,3,1,1,0,2,0,0,3,1).$$

\begin{lemma} \label{lemmata}
\begin{enumerate}
\item If $\alpha\in {\bf N}_{0}^{n+1}$ is a vector such that 
$\alpha_{n+1}\neq 0 $ and for all $1\leq
i\leq n+1$ either $\alpha_{i}=0$ or $\sum _{j=0}^{i}\alpha_{j}=i$,
then $\alpha=c(S)$ for only and only one subset $S\subseteq
\{1,\ldots,n \}$.
\item Any $\alpha$ as above is determined by the placement of its zeros.
\item Any $\alpha$ as above is determined by the value and order of
its nonzero elements.
\end{enumerate}
\end{lemma}

{\em Proof:} {\bf (1 and 2)} Let $S=\{i:\alpha_{i}=0 \}$.

{\bf (3)} Each nonzero element $a$ must be preceded by $a-1$ zeros.
Thus, the nonzero elements determine the zeros.$\Box $

We will write partitions with their parts arranged in lexicographical
order. (This is the usual way one write partitions, and is of
particular use, for example, in studying limited growth functions.)
The {\em type} $t(\phi )$ of a partition $\phi $ is the 
ordered sequence giving the 
length of each of its blocks when listed in lexicographical order.
For example, a partition $\phi $ of the integers from 1 to 20
according to 
their number of prime factors would be written like this: 
$$\phi =
\{\{1 \},\{2,3,5,7,11,13,17,19\}, \{4,6,9,10,14,15\}, \{8,12,18,20\},
\{16 \}\},$$    
and would thus be of type $t(\phi )=(1,8,6,4,1)$.   

For every set $S\subseteq \{1,\ldots,n \} $, we will associate a
collection of partitions $\Pi _{S}\subseteq \Pi _{n+1}$ as follows:
$\Pi _{S}$ 
consists of the set of all partitions $\phi \vdash \{1,\ldots,n+1 \}$
whose type when written backwards is the same as the code of $S$ when
written without its zeros.

\begin{lemma} \label{laugh}
 The set of classes $\{ \Pi _{S}: S\in {\cal B}_{n}\}$ forms a
partition of $\Pi _{n+1}$.   
\end{lemma}

{\em Proof:} Lemma \ref{lemmata}.$\Box $

\begin{prop}\label{18}
Let $S\rightarrow S'=S\cup \{i \}$ be a link in a de Bruijn's
decomposition of ${\cal
B}_{n}$. Then there is a simply described order-preserving  injection
of $\Pi _{S}$ into $\Pi _{S'}$. 
\end{prop}

{\em Proof:} We will explicitly construct the injection $f:\Pi
_{S}\rightarrow \Pi _{S'}.$ Part \ref{star} of Proposition
\ref{starprop} guarantees that 
difference between 
the codes $c(S)$ and $c(S')$ is that some substring of the form $k1$
is replaced by a string of the form $0(k+1)$. Let $\phi \in \Pi _{S}$.
Let $A=\{a_{1},\ldots,a_{k} \}$ and $B=\{b \}$ respectively be the
blocks of $\phi $ which correspond 
to the $k$ and the $1$ in the substring $k1$ which is being replaced.
Let $f(\phi ) $ be the replacement of $A$ and $B$ with their union:
$f(\phi ) = \phi \cup \{A\cup B \}- \{A,B \} $. Clearly, $f(\phi )$ is a
refinement of $\phi $, and $f(\phi )\in \Pi _{S'}$. Moreover, $f$ is
injective, since we can recover $\phi $ from $f(\phi )$. That is to
say, we can recover $A$ and $B$ from $A\cup B$, since
$b=\min (A\cup B).\Box $  

By applying Proposition \ref{18} to the entire decomposition, we
obtain the following theorem.

\begin{thm}
For each link $S\rightarrow S'$ in de Bruijn's decomposition of ${\cal
B}_{n}$ into symmetric chains, there are corresponding links between
all of $\Pi _{S}$ and part of $\Pi _{S'}$. Each resulting chain
beginning at rank $i$ ends at rank $n-i$ or higher. Thus, by pruning
the tops of these chains we achieve a partial decomposition of $\Pi
_{n+1}$ into symmetric chains. This decomposition is maximal since it
includes every partition with more than $\fl{(n+1)/2}$ parts.$\Box $
\end{thm}

The Stirling number $S_{nk}$ gives the number of partitions of an
$n$-set into $k$ blocks. 
\begin{cor}
The Stirling numbers of the second kind $S_{nk}$ obey the following
identities for all $n\geq 0$: 
$$S_{nn}\leq S_{n(n-1)}\leq \cdots \leq S_{n\fl{(n+1)/2}}$$
and for $k\leq n/2$
$$ S_{nk}\geq S_{n(n-k)}.\Box  $$
\end{cor}

\begin{example}
The de Bruijn decomposition of ${\cal B}_{3}$ consists of the chains
 $(\emptyset , \{1
\},\{1,2 \},\{1,2,3 \})$, $(\{2 \},\{2,3 \})$, and $(\{3 \},\{1,3
\})$.
Encoding the sets above yields the following partition types:
$(1111,112,13,4)$, $(121,31)$, and $(211,22)$. We then compute the
partitions of each type:
$$ \begin{array}{| r | r@{\rightarrow}l | lll |} \hline
S\in {\cal B}_{n}
	&	\multicolumn{2}{c|}{c(S)}  
			&	\multicolumn{3}{c|}{\Pi_{S}}\\ \hline \hline  
\emptyset &	1111 &	1111&	1/2/3/4&{}&	{}\\
\{1 \}&		0211 &	112 &	1/2/34&	{}&	{}\\
\{1,2 \} &	0031 &	13 &	1/234&	{}&	{}\\
\{1,2,3 \} &	0004 &	4 &	1234&	{}&	{}
                                                  \\[0.25in]
\{2 \} &	1021 &	121 &	1/23/4,&1/24/3&	{}\\
\{2,3 \} &	1003 &	31 &	123/4,&	124/3,&	134/2
                                                  \\[0.25in]
\{3 \} &	1102 &	211 &	12/3/4,&13/2/4,&14/2/3\\
\{1,3 \}&	0202 &	22 &	12/34,&	13/24,&	14/23\\ \hline
\end{array} $$
Reading down the columns marked $\Pi _{S}$ gives the chains in the
decomposition of $\Pi _{4}$.
\end{example}
 
\section{Other Applications of Set Coding}

The set coding defined above (Definition \ref{coding}) has been found
to have several applications. Indeed, in  \cite{DDN}, this coding was
used to give a simpler expression of the relationship between
elementary and complete homogeneous symmetric functions.

Let us write $e_{i}$ for the $i$th letter of $c(S)$ where
$S\subseteq \{1,2,\ldots,n-1 \}$. That is to
say, 
$$e_{i}=c(S)_{i}=\left\{ \begin{array}{ll}
0		&\mbox{if $i\in S$}\\
i-\sum _{j=1}^{i-1}c(S)_{j}		&\mbox{if $i\nin S$}
         \end{array}\right.$$
for $1\leq i\leq n$. Then we have \cite{DDN}
$$h_{n}=\sum _{\scriptstyle c(S)\atop \scriptstyle S\subseteq
\{1,2,\ldots,n-1 \}} (-1)^{|S|} a_{e_{1}}a_{e_{2}}\cdots a_{e_{n}}$$
where we denote by $h_{i}$ (resp. $a_{i}$) the $i$th complete (resp.
elementary) homogeneous symmetric function. 

Our coding also allows a straightforward derivation of the Bell
numbers' generating function. Let $g$ be a real function and let
$F(x)=e^{g(x)}. $ Now, let
$g^{(i)}(x)$ and $F^{(i)}(x)$ denote the $i$th derivatives of $F(x)$
and $g(x)$ 
respectively ($i\geq 0$). Then one can prove that
$$ F^{(n)}(x) = F(x) \left( 
\sum _{\scriptstyle c(S)\atop \scriptstyle S\subseteq \{1,2,\ldots,n-1 \} } 
\prod _{\scriptstyle i=1\atop \scriptstyle 
e_{i}\neq 0}^{n} {i-1\choose e_{i}-1 } g^{(i)}(x) \right). $$
(See \cite{DDN}.) 

Here we are interested in the case $g(x)=e^{x}-1$. In such an event,
$$ \exp(e^{x}-1) = \sum _{k=0}^{\infty } 
\sum _{\scriptstyle c(S)\atop \scriptstyle S\subseteq \{1,2,\ldots,n-1
\} } 
\prod _{\scriptstyle i=1\atop \scriptstyle e_{i}\neq 0}^{n} {i-1
\choose e_{i}-1 }\frac{x^{k}}{k!}.$$
Using lemmata \ref{lemmata} and \ref{laugh}, together with the concept
of the type $t(\phi )$ of a partition,  it can be shown that 
$$ B_{n}= \sum _{\scriptstyle c(S)\atop \scriptstyle S\subseteq
\{1,2,\ldots,n-1 \}} \prod _{\scriptstyle i=1\atop \scriptstyle
e_{i}\neq 0} {i-1\choose e_{i}-1}.$$

Note to compositor: Please be sure that the bold-face parentheses {\bf
()} are easily distinguished from regular parenthesis ().

 \end{document}